\documentclass[11pt]{article}
\usepackage{amscd,amsfonts,amsmath,amsthm}

\def\p{\partial}
\def\we{\wedge}
\def\half{{\textstyle\frac12}}
\def\tD{\tilde{D}}
\def\bW{\overline{W}}
\def\oW{\overline{W}}
\def\bomega{\overline{\omega}}
\def\eps{\varepsilon}
\def\RE{\mathbb{R}}
\def\CX{\mathbb{C}}
\def\CP{\mathbb{C}P}
\def\ZE{\mathbb{Z}}
\def\HE{\mathbb{H}}
\def\tphi{\tilde{\phi}}

\let\phi\varphi

\renewcommand{\Re}{\operatorname{Re}}

\newtheorem{thrm}{Theorem}
\newtheorem{prop}[thrm]{Proposition}
\theoremstyle{definition}
\newtheorem{defi}{Definition}
\newtheorem{example}{Example}
\theoremstyle{remark}
\newtheorem{cond}{Condition}

\newtheorem*{remarks}{Remarks}
\newtheorem*{acknowledge}{Acknowledgements}

\begin{document}

\title{Constructions of compact $G_2$-holonomy manifolds}
\author{{\sc Alexei Kovalev}\\[5pt]
DPMMS, University of Cambridge}
\date{}
\maketitle

\begin{abstract}
We explain the constructions for two geometrically different classes of
examples of compact Riemannian 7-manifolds with holonomy $G_2$. One method
uses resolutions of singularities of appropriately chosen 7-dimensional
orbifolds, with the help of asymptotically locally Euclidean spaces. Another
method uses the gluing of two asymptotically cylindrical pieces and requires a
certain matching condition for their cross-sections `at infinity'.
\end{abstract}

\section{Introduction}

The Lie group $G_2$ occurs as an exceptional case in Berger's
classification of the Riemannian holonomy groups, in
dimension 7. Riemannian manifolds with holonomy $G_2$ are Ricci-flat and
admit parallel spinor fields.
The purpose of these notes is to give an introduction
to two methods of producing examples of compact Riemannian
7-manifolds with holonomy group~$G_2$.

For a detailed introduction to $G_2$-structures on 7-manifolds and the
$G_2$ holonomy group we refer to \cite[Ch.~11]{salamon},
\cite[Ch.~10]{joyce-book} and the article by Karigiannis in this volume. Here
we briefly recall the foundational results that we need.

The Lie group $G_2$ may be defined as the stabilizer, in the action of
$GL(7,\RE)$, of the 3-form \cite[p.~539]{bryant}
\begin{equation}\label{phi0}
\phi_{0} = dx^{123} + dx^{145} + dx^{167} + dx^{246}
 - dx^{257} - dx^{347} - dx^{356} \in \Lambda^{3}(\RE^{7})^{*},
\end{equation}
where $x^k$ are the standard coordinates on $\RE^{7}$ and
$dx^{klm}=dx^k\we dx^l\we dx^m$. Every linear isomorphism of $\RE^7$
preserving $\phi_0$ also preserves the Euclidean metric
$\sum_{i=1}^7(dx^i)^2$ and orientation of $\RE^7$, thus $G_2$ is a subgroup
of~$SO(7)$. The $GL(7,\RE)$-orbit of $\phi_0$ is open in
$\Lambda^{3}(\RE^{7})^{*}$.

Let $M$ be a 7-dimensional manifold. Then every $G_2$-structure on $M$
is induced by a choice of a smooth differential 3-form $\phi$ such 
that for each $p\in M$ there is a linear isomorphism $\iota_p:\RE^7\to T_pM$
with $\iota_p^*(\phi(p))=\phi_{0}$. We say a 3-form $\phi$ is {\em positive}
when $\phi$ satisfies the latter condition and denote
by $\Omega^3_+(M)\subset\Omega^3(M)$ the subset of all
positive 3-forms on~$M$. Note that for a compact $M$ the subset
$\Omega^3_+(M)$ is open in the uniform norm topology.
We shall sometimes, slightly informally, say that a differential form
$\phi\in\Omega^3_+(M)$ {\em is} a $G_2$-structure on~$M$. 

We can see from the above that every $G_2$-structure $\phi\in\Omega^3_+(M)$
determines on~$M$ a metric $g(\phi)$ and an
orientation, hence also the Hodge star~${*}_\phi$.

\begin{thrm}[cf.~\cite{fernandez-gray}]
Let $M$ be a 7-manifold endowed with a $G_2$-structure
$\phi\in\Omega^3_+(M)$.
Then the following are equivalent.

(a) The holonomy of the metric $g(\phi)$ is contained in $G_2$.

(b) $\nabla\phi=0$, where $\nabla$ is the Levi--Civita connection of $g(\phi)$.

(c) \begin{equation}\label{torsion-free}
d\phi = 0,\qquad
d{*}_\phi \phi = 0,
\end{equation}

(d) The intrinsic torsion of the $G_2$-structure $\phi$ vanishes.
\end{thrm}

Note that the second equation in~\eqref{torsion-free} is non-linear because
${*}_\phi$ depends non-linearly on~$\phi$.

We say that $(M,\phi)$ is a {\em $G_2$-manifold} if $\phi$ is a positive
3-form satisfying~\eqref{torsion-free}. If, in addition, the holonomy of
$g(\phi)$ is all of~$G_2$, then we shall call $(M,\phi)$ an {\em irreducible
$G_2$-manifold}.
\begin{prop}[{\cite[Prop.~10.2.2]{joyce-book}}]\label{irreducible}
A compact $G_2$-manifold is irreducible if and only if $\pi_1(M)$ finite.
\end{prop}

A key idea in the known methods of constructing irreducible
$G_2$-manifolds is that one first achieves on $M$ a $G_2$-structure $\phi$
which is, in some sense, an `approximate' solution of~\eqref{torsion-free}
with $d\phi=0$ and 
$d{*}_\phi\phi$ having a small norm, in a suitable Banach space. In more
geometric terms, the $G_2$-structure $\phi$ then has small torsion. Then one
uses perturbative analysis to obtain a correction term $d\eta$, for a 2-form
$\eta$ small in the $C^1$ norm, so that $\phi+d\eta$ is a valid
$G_2$-structure and a solution of~\eqref{torsion-free}.

We shall explain methods of finding the desired approximate solutions
of~\eqref{torsion-free} by building compact Riemannian manifolds from `simpler 
pieces'. These will be non-compact or singular $G_2$-manifolds whose metrics
are flat or have holonomy $SU(2)$ or $SU(3)$, which are subgroups
of~$G_2$. These latter metrics can be obtained by using the Calabi--Yau
analysis or written explicitly. The manifolds are patched together in a
`compatible' way to achieve, on the resulting compact manifolds,
$G_2$-structures with arbitrarily small torsion.

More precisely, one obtains 1-dimensional families of metrics depending on
a certain `gluing parameter' taking values in a semi-closed interval. The
limits of these families may be interpreted as 
boundary points in a `partial compactification' of the $G_2$ moduli
space. (It is known that the moduli space of torsion-free
$G_2$-structures on a compact 7-manifold $M$ is a smooth manifold
of dimension the third Betti number $b^3(M)$.)

In these notes, we shall explain two ways of implementing the above strategy
with different respective limits in the boundary of the $G_2$ moduli space.

Recently, Joyce and Karigiannis \cite{JK} developed a new method of
constructing holonomy $G_2$ manifolds using analysis on families of
Eguchi-Hanson spaces. This construction is not reviewed here. It includes an
application of perturbative methods for $G_2$-structures with small torsion
but also requires significant additional methods to achieve a suitable small
torsion.

\begin{acknowledge}
These notes are an expanded version of lectures given in the Minischool on
$G_2$ manifolds at the Fields Institute, Toronto, in August 2017. I would like
to thank the organizers for inviting me.
\end{acknowledge}

\section{Construction by resolutions of singularities}

The method explained in this section was historically the first construction
of compact 7-manifolds with holonomy~$G_2$. It is due to Joyce
\cite{joyce-papers,joyce-book}.

Joyce's method produces one-parameter families of holonomy $G_2$
metrics $g_s$, $0<s\le\eps$. The limits of these families as $s\to 0$ can
be interpreted as boundary points in the $G_2$-moduli space and are given
by flat orbifolds. In particular, the limit spaces are {\em compact,
connected} and {\em singular}.

More precisely, the construction proceeds via the following steps.

\begin{enumerate}
\item
(a) Let $T^7=\RE^7/\ZE^7$ be the 7-torus with a flat $G_2$-structure
$\phi_0\in\Omega^3_+(T^7)$ induced from the standard
$G_2$-structure~\eqref{phi0} on the Euclidean $\RE^7$. Choose a finite
group $\Gamma$ of affine transformations of $\RE^7$ which preserve~$\phi_0$
and descend to diffeomorphisms of $T^7$. The quotient space
$(T^7/\Gamma)$ is an orbifold with a torsion-free $G_2$-structure,
still denoted by $\phi_0$, and a flat orbifold metric $g_0$ induced
by~$\phi_0$.
\\[5pt]
(b) For suitable choices of $\Gamma$, all the singularities of $T^7/\Gamma$
are locally modeled on $\RE^3\times(\CX^2/G)$ or $\RE\times(\CX^3/G)$, for
$G$ a finite subgroup of respectively $SU(2)$ or $SU(3)$, and can be resolved
using methods of complex algebraic geometry. Perform the resolutions to obtain
a smooth compact 7-manifold $M$ together with a resolution map
$\pi: M\to T^7/\Gamma$.
\item
(a) On $M$, one can `explicitly' define a 1-parameter family of closed
positive 3-forms $\phi_s\in\Omega^3_+(M)$, with $d\phi_s=0$ for $0<s\le\eps$,
such that the $G_2$-structures $\phi_s$ have small torsion. 
The forms $\phi_s$ converge as $s\to 0$ to $\pi^*\phi_0$ (respectively, the
induced metrics $g(\phi_s)$ converge to $\pi^* g_0$). One may also say that
the Riemannian manifolds $(M,g(\phi_s))$ converge in the Gromov--Hausdorff
sense to the flat orbifold $(T^7/\Gamma,g_0)$ as $s\to 0$.
\\[5pt]
(b) Apply perturbative analysis (more precisely, construct a convergent
sequence of iterations) to show that for every small $s>0$, the
$G_2$-structure $\phi_s$ can be deformed to a nearby torsion-free
$G_2$-structure $\tilde\phi_s$. If $\pi_1(M)$ is finite, then the holonomy
of the induced metric $\tilde{g}_s=g(\tilde\phi_s)$ is precisely the
group~$G_2$, i.e.\ $(M,\tilde\phi_s)$ is an irreducible $G_2$-manifold.
\end{enumerate}

We illustrate this method with an example taken from
\cite[\S 12.2]{joyce-book} (cf. also \cite{joyce-papers}) where some
technical details are relatively simple. Consider the group $\Gamma$
generated by
\begin{align*}
\alpha &: (x_{1}, \ldots, x_{7}) \mapsto
(\phantom{-}x_{1}, \phantom{-}x_{2}, \phantom{-}x_{3}, -x_{4}, 
\phantom{\half}{-}x_{5}, \phantom{\half}{-}x_{6}, \phantom{\half}{-}x_{7}) , \\
\beta &: (x_{1}, \ldots, x_{7}) \mapsto
(\phantom{-}x_{1}, -x_{2}, -x_{3}, \phantom{-} x_{4},
\phantom{\half{-}} x_{5}, \half {-} x_{6}, \phantom{\half}{-}x_{7}) , \\
\gamma &: (x_{1}, \ldots, x_{7}) \mapsto
(-x_{1}, \phantom{-}x_{2}, -x_{3}, \phantom{-} x_{4},
\half{-}x_{5}, \phantom{\half{-}} x_{6}, \half{-}x_{7}).
\end{align*}
The maps $\alpha,\beta,\gamma$ commute and each has order~2, thus $\Gamma$ is
isomorphic to~$\ZE_2^3$. The elements of $\Gamma$ descend to $T^7$ and
preserve $\phi_0$, making the quotient $T^7/\Gamma$ into a $G_2$-orbifold.

One can further check that the only elements of $\Gamma$ having fixed
points are $\alpha,\beta,\gamma$, each fixes 16 copies of~$T^3$ and these
are all disjoint. The subgroup generated by $\beta,\gamma$ acts freely on
the 16 tori fixed by~$\alpha$, so these correspond to 4 copies of $T^3$ in
the singular locus of $T^7/\Gamma$. Similar properties hold for the tori
fixed by $\beta$ and by~$\gamma$. Thus the singular locus $S$ of
$T^7/\Gamma$ is 12 disjoint copies of~$T^3$.  A neighbourhood of each
3-torus component of~$S$ is diffeomorphic to $T^3\times (\CX^2/(\pm 1))$.

The blow-up $\sigma:Y\to\CX^2/(\pm 1)$ at the origin resolves the singularity
with a complex surface $Y$ biholomorphic to $T^*\CP^1$, with the
exceptional divisor $E=\sigma^{-1}(0)\cong\CP^1$ corresponding to the zero 
section of $T^*\CP^1$. The canonical bundle of $Y$ is trivial and there is
a family of Ricci-flat K\"ahler
metrics $h_s$ on~$Y$ with holonomy equal to $SU(2)$ depending on a real
parameter $s>0$. The K\"ahler form of the metric $h_s$ may be written as
$\omega_s=\sigma^*(i\p\bar\p f_s)$, where
$$
f_s=\sqrt{r^4+s^4}+2s^2\log r - s^2\log(\sqrt{r^4+s^4}+s^2),
$$
$r^2=z_1\bar{z}_1+z_2\bar{z}_2$ and $(z_1,z_2)\in\CX^2$. The radius
function $r$ makes sense as a smooth function on $Y\setminus E$ and the
values of this function near $E$ can be interpreted as the distance to $E$
in the metric $h_s$. The forms $\omega_s$ extend smoothly over the exceptional
divisor~$E\subset Y$, thus the metrics $h_s$ are well-defined on $Y$. These
are the well-known Eguchi--Hanson metrics~\cite{EH}.

Comparing, for each $s>0$, the K\"ahler potential $f_s$ of $h_s$ with the 
K\"ahler potential $r^2$ of the Euclidean metric $h_0$ on $\CX^2$
we see that 
\begin{equation}\label{ALE2}
\nabla^k(h_s-h_0)=O(r^{-4-k})\quad\text{as } r\to\infty,
\quad\text{for all }k=0,1,2,\ldots\; ,
\end{equation}
which means that $h_s$ is an asymptotically locally Euclidean (ALE)
metric on~$Y$.

For each $\lambda>0$, the dilation map $Y\to Y$ induced by $(z_1,z_2)\mapsto
\lambda (z_1,z_2)$ pulls back $\omega_s$ to $\lambda^2\omega_{\lambda s}$.
It follows that $s$ is proportional to the diameter of the exceptional
divisor. One can further check that the injectivity radius of the
Eguchi--Hanson metric $h_s$ is proportional to~$s$ and that the uniform
norm of the Riemannian curvature is proportional to~$s^{-2}$.

Every Ricci-flat K\"ahler metric $h$ on a complex surface is in fact
hyper-K\"ahler: in addition to the original complex structure $I$ there are
(integrable) complex structures $J$ and $K$ satisfying quaternionic relations
$IJ=-JI=K$. For each $p\in Y$, there is an $\RE$-linear isomorphism
$\RE^4\to T_pY$ such that the linear maps $I(p),J(p),K(p)$ correspond to
multiplication by the unit quaternions $i,j,k$ via the standard
identification $\RE^4\cong\HE=\langle 1,i,j,k\rangle$ with the algebra of
quaternions. Also, the metric $h$ is K\"ahler with respect to each
$I,J,K$. We shall denote by $\kappa_I$, $\kappa_J$, $\kappa_K$ the
respective K\"ahler forms.

For a 3-torus $T^3$ with coordinates $x_1,x_2,x_3$, with a flat
metric $dx_1^2+dx_2^2+dx_3^2$ and a hyper-K\"ahler 4-manifold $Y$ as above,
the Riemannian product $T^3\times Y$ has holonomy in $SU(2)$. The product
metric is induced by a torsion-free $G_2$ structure on $T^3\times Y$, which is
\begin{equation}\label{su2}
\phi_{SU(2)}=dx_1\we dx_2\we dx_3 +
dx_1\we\kappa_I + dx_2\we\kappa_J - dx_3\we\kappa_K.
\end{equation}

We now define, for every small $\eps>0$, a smooth compact 7-manifold
$M=M_\eps$ by replacing a neighbourhood $T^3\times\{r<2\eps\}$ of each
3-torus component in the singular locus of $T^7/\Gamma$ by $T^3\times U$,
where $U=\sigma^{-1}(r<2\eps)\subset Y$ is a neighbourhood of the
exceptional divisor in~$Y$.
(Note that the manifolds $M_\eps$ are diffeomorphic to each other.)

On each $T^3\times U$ in~$M$, we smoothly interpolate, for $\eps<r<2\eps$,
between the flat $G_2$-structure $\phi_0$ induced from $T^7/\Gamma$ on the
complement of the regions $T^3\times U$ and the product $G_2$-structure
arising as in~\eqref{su2} from the appropriately rescaled
Eguchi--Hanson hyper-K\"ahler $h_s$ on $\sigma^{-1}(r<\eps)\subset U$. 
The ALE property of the Eguchi--Hanson metric allows to take the product
$G_2$-structure on $T^3\times Y$ to be asymptotic to the flat
$G_2$-structure on $T^3\times (\CX^2/(\pm 1))$.
We can obtain, for each sufficiently small $s>0$, a well-defined
positive 3-form $\phi_s$ on $M$ noting also that $\Omega^3_+(M)$
is an open subset of 3-forms in the uniform norm.
Furthermore, we can choose these $G_2$ \mbox{3-forms} on $M$ to be
{\em closed}, $d\phi_s=0$. 
Thus the $G_2$-structure $\phi_s$ is torsion-free away from the interpolation
region $\{\eps<r<2\eps\}$ but $\phi_s$ is not co-closed in that region.

The positive 3-forms $\phi_s$ are intended as `approximate solutions' of
the torsion-free equations~\eqref{torsion-free}, as $s\to 0$. The
parameter $s$ may be interpreted geometrically as the maximal diameter of
the pre-image of a singular point in $T^7/\Gamma$ under the resolution map
$M\to T^7/\Gamma$. We would like to perturb $\phi_s$ to actual solutions on
$M$. To this end, the following two conditions satisfied by $\phi_s$ are
important, cf.~\cite[Thm. 11.5.7]{joyce-book}.

\begin{cond}
One can construct a smooth 3-form $\psi_s$ on $M$
such that $d^*\phi_s=d^*\psi_s$ and
\begin{equation}\label{estim}
\|\psi_s\|_{L^2}<A_1s^4,  \|\psi_s\|_{C^0}<A_1s^3,
\text{ and } \|d^*\psi_s\|_{L^{14}}<A_1s^{16/7}.
\end{equation}
\end{cond}

\begin{cond}
The injectivity radius $\delta(g_s)$ and the Riemann curvature $R(g_s)$ of
the metric $g_s=g(\phi_s)$ on~$M$ satisfy the estimates
\begin{equation}\label{metric}
\delta(g_s)>A_2\,s,
\qquad
\|R(g_s)\|_{C^0}<A_3s^{-2}.
\end{equation}
\end{cond}

The construction of $\psi_s$ exploits the asymptotic and scaling properties
of the $G_2$-structure~\eqref{su2} on $T^3\times U$ `approximating' the flat
$G_2$-structure on $T^3\times\CX^2/(\pm 1)$.
The estimates~\eqref{metric} follow from the properties of the metric $h_s$
around the exceptional divisor on $Y$, which give the dominant
contributions for small~$s$.
In the conditions (i) and (ii) the norms and the formal adjoint $d^*$ are
taken with respect to the metric $g_s=g(\phi_s)$. The constants
$A_1,\,A_2,\,A_3$ are independent of $s$.

We can now state the existence result for torsion-free $G_2$-structures.
\begin{thrm}[{cf.\ \cite[Thm.~11.6.1]{joyce-book}}]\label{g2-exist}
Let $M$ be a compact 7-dimensional manifold and $\phi_s\in\Omega^3_+(M)$,
$0<s\le s_0$, a family of $G_2$-structures such that $d\phi_s=0$ and the
conditions (i) and (ii) above hold for all $s$.

Then there is an $\eps_0>0$ so that for each $s$ with $0<s\le\eps_0$ the
manifold $M$ admits a torsion-free $G_2$ structure
$\tilde\phi_s\in\Omega^3_+(M)$ in the same cohomology class as
$\phi_s$ and satisfying $\|\tilde\phi_s-\phi_s\|_{C^0}<Ks^{1/2}$ with some
constant $K$ independent of~$s$.
\end{thrm}

We next outline the proof of Theorem~\ref{g2-exist} following
\cite[p.236--237]{joyce-book2}, dropping the subscripts $s$ to ease the
notation. The desired torsion-free $G_2$-structure $\tphi=\tilde\phi_s$
will be obtained in the form
$\tphi=\phi+d\eta$, where $d\eta$ has a small uniform norm, so $\tphi$ is
a closed positive 3-form. We then need to satisfy the co-closed condition
$d^*_{\tphi}\tphi=0$ and this amounts to solving for a 2-form $\eta$ a
non-linear elliptic PDE which may be written as
\begin{equation}\label{eta}
d^*d\eta=-d^*\psi+d^*F(d\eta)
\end{equation}
where $F$ satisfies a quadratic estimate. A solution of~\eqref{eta}
is achieved by using iterations to construct a sequence
$\{\eta_j\}_{j=0}^\infty$ with $\eta_0=0$ and
$$
d^*d\eta_{j+1}=-d^*\psi+d^*F(d\eta_j), \qquad d^*\eta_{j+1}=0.
$$
One first argues that the sequence $\eta_j$ converges.

The proof of convergence is based on the following inductive estimates (all
the constants $C_i$ below are independent of $s$)
\begin{subequations}\label{induct}
\begin{gather}
\|d\eta_{j+1}\|_{L^2} \le
\|\psi\|_{L^2}+C_1\|d\eta_{j}\|_{L^2}\|d\eta_{j}\|_{C^0},
\\
\|\nabla d\eta_{j+1}\|_{L^{14}} \le C_2(\|d^*\psi\|_{L^{14}}+\|\nabla
d\eta_{j}\|_{L^{14}}\|d\eta_{j}\|_{C^0}+s^{-4}\|d\eta_{j+1}\|_{L^2}),
\\
\|d\eta_{j}\|_{C^0} \le
C_3(s^{1/2}\|\nabla d\eta_j\|_{L^{14}}+s^{-7/2}\|d\eta_j\|_{L^2}).
\end{gather}
\end{subequations}
The estimate~(\ref{induct}a) is proved by taking the $L^2$ product of both
sides with $\eta_{j+1}$ and integrating by parts, noting also condition~(i)
above. The proof of~(\ref{induct}b) uses an elliptic regularity estimate
for the operator $d+d^*$ considered for 3-forms on small balls on $M$ with
radius of order~$s$. The condition (ii) is also required here and
in~(\ref{induct}c) which uses the Sobolev embedding of $L^{14}_1$ in $C^0$
in dimension~7 and is again achieved by working on small balls with radius
of order~$s$.

For every sufficiently small $s$, we deduce from~\eqref{induct} 
that if $d\eta_j$ satisfies
\begin{equation}\label{bounds}
\|d\eta_j\|_{L^2}\le C_4s^4,\quad
\|\nabla d\eta_{j}\|_{L^{14}} \le C_5,\quad
\|d\eta_j\|_{C^0} \le Ks^{1/2},
\end{equation}
then these latter estimates hold for $d\eta_{j+1}$ and, by induction, for
all~$j$. Thus $d\eta_j$ is a bounded sequence in the $L^{14}_1$ norm
on~$\Lambda^3 T^* M$ and one can further show that $d\eta_j$ is a
Cauchy sequence. Further, we are free to assume that the forms $\eta_j$ are
in the $L^2$-orthogonal complement $\mathcal{H}^\bot$ of harmonic forms. As the
elliptic operator $d+d^*$ is bounded below on $\mathcal{H}^\bot$ it follows
that the sequence $\eta_j$ converges in the $L^{14}_2$ norm. In particular,
the last inequality of~\eqref{bounds} holds for the limit~$\eta$.

Finally, a careful elliptic regularity argument shows that $\eta$ is in fact
a smooth solution of~\eqref{eta}, thus completing the proof of
Theorem~\ref{g2-exist}.

The metrics on $M$ induced by $\phi_s$ have holonomy in~$G_2$ and it
remains to verify that the holonomy does not reduce further to a
subgroup of~$G_2$.
In the present case, the orbifold $T^7/\Gamma$ is simply-connected,
therefore $M$ is so, by the properties of the blow-up.
Thus $(M,\phi_s)$ is an irreducible $G_2$-manifold by
Proposition~\ref{irreducible}.

The discussed example may be considered as a generalization of the Kummer
construction of hyper-K\"ahler metrics of holonomy $SU(2)$ on K3 surfaces.

It is convenient to obtain the Betti numbers of $M$; these are determined
by $b^2(M)$ and $b^3(M)$. By considering the $\Gamma$-invariant classes in
$H^*_{\mathrm{dR}}(T^7/\Gamma)$ we obtain $b^2(T^7/\Gamma)=0$ and
$b^3(T^7/\Gamma)=7$. When resolving the singularities, we replaced a
deformation retract of $T^3$ with $T^3\times Y$ which is homotopy
equivalent to $T^3\times\CP^1$. Let $S$ denote the singular locus of
$T^7/\Gamma$. Comparing the cohomology long exact
sequence for the pairs $(T^7/\Gamma,S)$ and
$(M,\sqcup_{i=1}^{12}(T^3\times U))$, we find that each of the 12 instances 
of a resolution adds $b^i(T^3\times Y)-b^i(T^3)$ to the $i$-th Betti number
of~$M$. Thus $b^2(M)=12\cdot 1$ and $b^3(M)=7+12\cdot 3=43$.

Further examples of irreducible $G_2$-manifolds arise by using the above
method with different choices of finite groups $\Gamma$ and different
choices of resolutions of singularities of $T^7/\Gamma$. If every component
of the singular locus of $T^7/\Gamma$ has a neighbourhood diffeomorphic to
$T^3\times(\CX^2/G)$ for a finite subgroup $G$ of $SU(2)$ or to
$S^1\times(\CX^3/G)$ for a finite subgroup $G$ of $SU(3)$ acting freely on
$\CX^3\setminus\{0\}$, then it is known from complex algebraic geometry
that one can find crepant resolutions, $\sigma_2:Y_2\to\CX^2/G$ or
$\sigma_3:Y_3\to\CX^3/G$ respectively, with the canonical bundle of $Y_i$
holomorphically trivial.

The Ricci-flat K\"ahler (thus hyper-K\"ahler) metrics on the complex
surfaces $Y_2$ asymptotic to $\CX^2/G$ in the sense of~\eqref{ALE2}, for
each $G$, were constructed by Kronheimer~\cite{kronh} using hyper-K\"ahler
quotients.

In complex dimension~3, the existence of ALE Ricci-flat holonomy $SU(3)$
metrics on $Y_3$ asymptotic to $\CX^3/G$ follows from the solution
of ALE version of the Calabi conjecture, see \cite[Ch. 8]{joyce-book}
and references therein.
The asymptotic rate for the metrics $h$ is given by
$$
\nabla^k(h-h_0)=O(r^{-6-k})\quad\text{as } r\to\infty,
\quad\text{for all }k=0,1,2,\ldots\; ,
$$
where $h_0$ is the pull-back of the Euclidean metric on $\CX^3/G$.
The K\"ahler forms of $h$ and $h_0$ satisfy
$$
\omega-\omega_0=i\p\bar\p u,\quad
\nabla^k u = O(r^{-4-k})\quad\text{as } r\to\infty
$$
(cf.~\cite[Thm. 8.2.3]{joyce-book}).
The holonomy being $SU(3)$ means there is a choice of nowhere vanishing 
(3,0)-form $\Omega$ on $Y_3$ (sometimes called a holomorphic volume form),
such that $\omega^3/3!=(i/2)^3\Omega\we\bar\Omega$.
A torsion-free $G_2$-structure on $S^1\times Y_3$ defined by
\begin{equation}\label{su3}
\phi_{SU(3)}=dx\we\omega+\Re\Omega,
\end{equation}
induces a product metric corresponding to $(dx)^2$ and $\omega$,
where $x$ is the usual coordinate on $S^1=\RE/\ZE$.

The singularities of $T^7/\Gamma$ can be resolved with copies of
$T^3\times U_2$ or $S^1\times U_3$ (where $U_i$ is a neighbourhood of
$\sigma_i^{-1}(0)$ in $Y_i$) in a manner similar to the example above.
One obtains compact smooth 7-manifolds $M$ and closed positive 3-forms
$\phi_s$ on~$M$ satisfying the hypotheses of Theorem~\ref{g2-exist}.
More generally, the method extends to situations when the singularities
of $T^7/\Gamma$ are only {\em locally} modeled on $\RE^3\times(\CX^2/G)$
or $\RE\times(\CX^3/G)$. In the latter case, $G$ need not act freely on
$\CX^3\setminus\{0\}$ resulting in a more complicated singular locus of
$T^7/\Gamma$.

Joyce found a large number of orbifolds $T^7/\Gamma$ with suitable
resolutions of singularities. In particular, 252 examples of topologically
distinct compact 7-manifolds admitting holonomy $G_2$ metrics are worked
out in \cite[Ch. 11]{joyce-book}, including some manifolds with non-trivial
fundamental group. The Betti numbers of these examples are in the range
$0\le b^2\le 28$ and $4\le b^3\le 215$. There is a evidence that many more
further topological types can be constructed by the same method.

\section{Construction by generalized connected sums}

The method of constructing compact holonomy $G_2$ manifolds discussed in
this section is sometimes called a `twisted connected sum'. The construction
was originally developed by the author in \cite{g2paper} and included an
important idea due to Donaldson. Generalizations and many new examples
appeared in \cite{k3g2,CHNP1,CHNP2,extra-twisted}.

The connected sum construction produces one-parameter families of holonomy
$G_2$ metrics $g_t$, $T_0\le T<\infty$, on compact manifolds with `long
necks'. The parameter $t$ here is asymptotic, as $T\to\infty$, to the diameter
of the metric $g_t$. We may think of the respective families of torsion-free
$G_2$-structures as paths in the $G_2$ moduli space, going to the boundary
as one `stretches the neck', the limit boundary point corresponding to the
disjoint union of the initial two asymptotically cylindrical pieces. So, in
this construction, the limit spaces are {\em disconnected, non-compact} and
{\em smooth}.

A twisted connected sum is an instance of generalized connected sum of a
pair of asymptotically cylindrical Riemannian manifolds which, in the
present case, are $G_2$-manifolds. The asymptotically cylindrical
$G_2$-manifolds we require are Riemannian products $W\times S^1$, where $W$
is a Ricci-flat K\"ahler manifold 
with cylindrical end asymptotic to a Riemannian product $D\times S^1\times
[0,\infty)$ with $D$ a K3 surface with a hyper-K\"ahler metric.
For certain pairs of the K3 surfaces $D_1,\; D_2$ there is a way to `join'
the two latter asymptotically cylindrical manifolds at their ends.
We obtain a compact simply-connected manifold $M$ and a $G_2$-structure with
small torsion on~$M$ to which a perturbative analysis can be applied.

We now describe the key steps in the construction in more detail, starting
with the asymptotically cylindrical Calabi--Yau threefolds~$W$.

\begin{thrm}[\cite{TY}, \cite{g2paper}, \cite{HHN}]\label{bCY}
Let $\bW$ be a compact K\"ahler threefold with K\"ahler form $\bomega$
and suppose that a K3 surface $D\in |-K_{\bW}|$ is an anticanonical divisor
on $\bW$ with holomorphically trivial normal bundle $N_{D/\bW}$. Denote by
$z$ a complex coordinate around~$D$ vanishing to order one precisely on~$D$.
Suppose that $\bW$ is simply-connected and the fundamental group of
$W=\bW\setminus D$ is finite. 

Then $W$ admits a complete Ricci-flat K\"ahler metric, with holonomy~$SU(3)$,
with K\"ahler form $\omega$ and a non-vanishing holomorphic
$(3,0)$-form~$\Omega$. These are asymptotic to the product
cylindrical Ricci-flat K\"ahler structure on $D\times S^1\times \RE_{>0}$
\begin{gather*}
\omega=\kappa_I+dt\we d\theta+d\psi,
\\
\Omega=(\kappa_J+i\,\kappa_K)\we (dt+id\theta)+d\Psi,
\end{gather*}
where $\exp(-t-i\theta)=z$, for $(\theta,t)\in S^1\times \RE_{>0}$ and the
forms $\psi,\Psi$ exponentially decay as $t\to\infty$.
Also $\kappa_I$ is the Ricci-flat K\"ahler metric
on~$D$ in the class $[\bomega|_D]$ and $\kappa_J+i\kappa_K$ is a
non-vanishing holomorphic $(2,0)$-form on~$D$.
\end{thrm}
\begin{remarks}
Any threefold $\bW$ satisfying the hypotheses of Theorem~\ref{bCY} is
necessarily projective and algebraic~\cite[Prop. 2.2]{k3g2}. The holomorphic
coordinate $z$ extends to a meromorphic function $\bW\to\CP^1$ vanishing
precisely on $D$.

Theorem~\ref{bCY} extends to higher dimensions $m\ge 3$ with $D$ replaced by
a compact simply-connected Calabi-Yau $(m-1)$-fold. The result may be
regarded as a solution of an `asymptotically cylindrical version' of the
Calabi conjecture.
\end{remarks}

It will be convenient to extend the parameter $t$ along the cylindrical end in
Theorem~\ref{bCY} to a smooth function $t$ defined on all of~$W$ with
$t<0$ away from a tubular neighbourhood of~$D$.
We shall also assume that the holomorphic 2-form on a K\"ahler K3 surface $D$
is normalized so that $\kappa_I^2=\kappa_J^2=\kappa_K^2$, with the implied
normalization of a holomorphic 3-form $\Omega$ on~$W$. The Ricci-flat K\"ahler
(hyper-K\"ahler) structure on $D$ is in fact determined by the triple
$\kappa_I,\kappa_J,\kappa_K$ (cf.~\cite[p. 91]{hitchin}.

The following relation between K3 surfaces is crucial for the connected sum
construction of $G_2$-manifolds.
\begin{defi}\label{matching}
We say that two Ricci-flat K\"ahler K3 surfaces 
$(D_1,\kappa'_I,\kappa'_J+i\kappa'_K)$,
$(D_2,\kappa''_I,\kappa''_J+i\kappa''_K)$
satisfy the {\em Donaldson matching condition} if there exists an isometry
of lattices $h:H^2(D_2,\ZE)\to H^2(D_1,\ZE)$, so that the $\RE$-linear
extension of $h$ satisfies
\begin{equation}\label{rotation}
h:\; [\kappa''_I]\mapsto [\kappa'_J],\qquad
[\kappa''_J]\mapsto [\kappa'_I],\qquad
[\kappa''_K]\mapsto [-\kappa'_K].
\end{equation}
\end{defi}
It follows, by application of the Torelli theorem for K3 surfaces,
that there is a smooth map
$$
f:D_1\to D_2, \text{ such that }h=f^*.
$$
Note that $f$ is {\em not} a holomorphic map between $D_1$ and $D_2$ (with
complex structures~$I$), though $f$ is an isometry of
the underlying Riemannian 4-manifolds. In particular, the pull back $f^*$
rotates the 2-forms of the hyper-K\"ahler triple (not just their cohomology
classes),
$\kappa''_I\mapsto\kappa'_J$, $\kappa''_J\mapsto\kappa'_I$,
$\kappa''_K\mapsto -\kappa'_K$.

Now if $(W,\omega,\Omega)$ is an asymptotically cylindrical Calabi-Yau
manifold given by Theorem~\ref{bCY}, then $W\times S^1$ 
has a torsion-free $G_2$-structure given by~\eqref{su3}
$$
\phi_{W}=d\tilde\theta\we\omega+\Re\Omega,
$$
where $\tilde\theta$ is the standard coordinate on the $S^1$ factor.
The form $\phi_{W}$ is asymptotic to a cylindrical product torsion-free 
$G_2$-structure $\phi_{\infty}$ on the cylindrical end
$D\times [0,\infty)\times S^1\times S^1\subset W\times S^1$,
$$
\phi_{\infty}=dt\we d\theta\we d\tilde\theta +
d\tilde\theta\we\kappa_I + dt\we\kappa_J - d\theta\we\kappa_K.
$$
corresponding to the hyper-K\"ahler structure
$(\kappa_I,\kappa_J,\kappa_K)$ on $D$ (cf. \eqref{su2}).

For $i=1,2$ and $T>0$, let $W_{i,T}$ be a compact manifold with boundary
obtained by truncating $W_i$ at $t_i=T+1$ (where $t_i$ is the parameter
along the cylindrical end as in Theorem~\ref{bCY}). We can smoothly cut off
each $\phi_{W_i}$ to obtain on $W_{i,T}$ a closed $G_2$-structure $\phi_{W,T}$
so that $\phi_{W_i,T}$ equals its cylindrical asymptotic model $\phi_\infty$ on
a collar neighbourhood $D_i\times S^1\times S^1\times [T,T+1]$ of the
boundary. 

Suppose that $D_1$ and $D_2$ satisfy the Donaldson matching condition.
Then we can define a compact 7-manifold
\begin{equation}\label{7mfd}
M=M_T=(W_{1,T+1}\times S^1)\cup_F (W_{2,T+1}\times S^1)
\end{equation}
by identifying the collar neighbourhoods of the boundaries using a map
\begin{equation}\label{neck}
\begin{split}
F:D_1\times S^1\times S^1\times [T,T+1]&\to
D_2\times S^1\times S^1\times [T,T+1],
\\
(y,\theta,\tilde\theta,T+t)&\mapsto (f(y),\tilde\theta,\theta,T+1-t).
\end{split}
\end{equation}
The form $\phi_{\infty}|_{[T,T+1]}$ is preserved by $F$, so the
$G_2$-structures $\phi_{i,T}$ agree on the overlap and patch together to a 
well-defined closed 3-form $\phi_T$ on~$M$. It is easy to see that
$\phi_T$ is a well-defined $G_2$-structure on~$M$ for every large~$T$.

Another important property of the map $F$ is that $F$ identifies the $S^1$
factor in $W_{1,T+1}\times S^1$ with a circle around the divisor in the other
threefold $W_{2}$ and {\it vice versa}. This eliminates the possibility of an
infinite fundamental group of~$M$, in particular, $M$ will be simply-connected
when the threefolds $W_1$ and $W_2$ are so.

The $G_2$-structure form on $M$ satisfies $d\phi_T=0$,
one of the two equations in~\eqref{torsion-free}, but the co-derivative
$d*_T\phi_T$ in general will not vanish. The cut-off functions introduce
`error terms' which depend on the difference between the $SU(3)$-structures
on the end of $W_i$ and on its cylindrical asymptotic model, and can be
estimated as
$$
\|d*_T\phi_T\|_{L^p_k} < C_{p,k} e^{-\lambda T},
$$
with $\lambda>0$. Here $*_T$ denotes the Hodge star of the metric $g(\phi_T)$.

The next result shows that for a sufficiently long neck the $G_2$-structure
$\phi_T$ on $M$ can be made torsion-free by adding a small correction term.
	\begin{thrm}\label{gluing}
Suppose that each of $\oW_1,D_1$ and $\oW_2,D_2$ satisfies the hypotheses
of Theorem~\ref{bCY} and the $K3$ surfaces $D_j\in|-K_{\oW_j}|$ satisfy the
Donaldson matching condition. Let $M$ be the compact 7-manifold $M$ defined
in~\eqref{7mfd} with a closed $G_2$-structure $\phi_T$ induced from
$\phi_{W_1},\, \phi_{W_2}$.

Then $M$ has finite fundamental group.
Furthermore, there exists $T_0\in\RE$ and for every $T\ge T_0$ a unique
smooth 2-form $\eta_T$ on~$M$, orthogonal to the closed forms, so that the
following holds.
\\
(a) $\|\eta_T\|_{C^1}< A\cdot e^{-\mu T}$, for some constants
$A,\mu>0$ independent of~$T$, where the norm is defined using the
metric $g(\phi_T)$. In particular, $\phi_T+d\eta_T$ is a valid
$G_2$-structure on $M$.
\\
(b) The closed 3-form $\phi_T+d\eta_T$ satisfies
	\begin{equation}\label{2nd}
d*_{\phi_T+d\eta_T}(\phi_T+d\eta_T)=0.
	\end{equation}
and so $\phi_T+d\eta_T$ defines a metric with holonomy $G_2$ on~$M$.
	\end{thrm}

As discussed in the previous section, the perturbative problem~\eqref{2nd}
can be rewritten as a \mbox{non-linear} elliptic PDE for the 2-form~$\eta$.
When $\eta$ has a small norm this PDE takes the form 
$a(\eta)=a_0+A\eta+Q(\eta)=0$, where $a_0=d*_T\phi_T$, the linear
elliptic operator $A=A_T$ is a compact perturbation of the Hodge Laplacian
of the form $dd^*+d^*d+O(e^{-\eps T})$, $\eps>0$ and $Q(\eta)$ satisfies a
quadratic estimate in $d\eta$.

One can use elliptic theory for manifolds with cylindrical ends and the
gluing analysis for the problem at hand is then simplified, compared to the
general situation of Theorem~\ref{g2-exist}.
The central idea in the proof of Theorem~\ref{gluing} may be informally
described as follows. For small $\eta$, the map $a(\eta)$ is approximated
by its linearization and so there would be a unique small solution $\eta$ to
the equation $a(\eta)=0$, for every small $a_0$ in the range of~$A$.
This perturbative approach requires the invertibility of $A$ and a
suitable upper bound on the operator norm $\|A_T^{-1}\|$, as
$T\to\infty$. This bound determines what is meant by `small~$a_0$' above.

As we actually need the value of $d\eta$ rather than $\eta$ we may consider
the equation for $\eta$ in the orthogonal complement of harmonic 2-forms
on~$M$ where the Laplacian is invertible.
We use the technique similar to~\cite[\S 4.1]{ks} based on Fredholm theory for
the asymptotically cylindrical manifolds and weighted Sobolev spaces
to find an upper bound $\|A_T^{-1}\|<Ge^{\delta T}$. Here the
constant $G$ is independent of $T$ and $\delta>0$ can be taken
arbitrarily small. So, for large~$T$, the growth of $\|A_T^{-1}\|$ is
negligible compared to the decay of $\|d*_T\phi_T\|$ and the `inverse
function theorem' strategy applies to give the required small
solution~$\eta_T$ in a (appropriately chosen) Sobolev space. Standard elliptic
methods show that this~$\eta_T$ is in fact smooth.

\subsection{Some examples and further results}

In order to make irreducible $G_2$-manifolds using the connected sum
construction, we require pairs $\bW_1$, $\bW_2$ of complex
algebraic threefolds with matching anticanonical K3 divisors
$D_i\subset\bW_i$. We begin with an example based on some classical algebraic
geometry.

\begin{example}\label{one}
The intersection of three generically chosen quadric hypersurfaces in
$\CP^6$ defines a smooth K\"ahler threefold $X_8$. It is
simply-connected and the characteristic class $c_1(X_8)$ of its
anticanonical bundle is the pull-back to~$X_8$ of the positive
generator of the cohomology ring $H^*(\CP^6)$. This tells us that the
anticanonical bundle $K^{-1}_{X_8}$ is the restriction to $X_8$ of the
tautological line bundle $\mathcal{O}(1)$ over $\CP^6$.  It follows
that any anticanonical divisor $D$ on $X_8$ is obtained by taking an
intersection $D=X_8\cap H$ with a hyperplane $H$ in $\CP^6$. A generic
such hyperplane section $D$ is a complex surface, isomorphic to
a smooth complete intersection of three quadrics in $\CP^5$. This
is a well-known example of a K3 surface.

Conversely, starting from a smooth intersection $D$ of three quadrics in
$\CP^5$ we can write down a smooth threefold $X_8\subset\CP^6$ as above
containing the K3 surface $D$ as an anticanonical divisor. 

Consider another anticanonical divisor $D'=X_8\cap H'$ and let
$\tilde{X}_8\to X_8$ be the blow-up of the second hyperplane section
$C=D\cap D'=X_8\cap H\cap H'$. (It is convenient, though not strictly
necessary, to choose $D'$ so that $C$ is a non-singular connected complex
curve.) The pencil defined by $D$ and $D'$ lifts, via the proper transform, to
a pencil consisting of the fibres of a holomorphic map
$\tilde{X}_8\to\CP^1$. In particular, the 
K3 divisor $D$ lifts to an isomorphic K3 surface $\tD$ which is an
anticanonical divisor in $\tilde{X}_8$ and has trivial normal bundle.
Moreover, a K\"ahler metric on $\tilde{X}_8$ may be chosen so that
$\tilde{D}$ and~$D$ are isometric K\"ahler K3 surfaces.

It is not difficult to check that $\tilde{X}_8\setminus\tD$ is
simply-connected, noting that $\tD$ and $X_8$ are so and considering an
exceptional curve in the blow up $\tilde{X}_8$. The pair
$\tilde{X}_8$, $\tilde{D}$ thus satisfies all the hypotheses of
Theorem~\ref{bCY}, and so the quasiprojective threefold
$W=\tilde{X}_8\setminus\tilde{D}$ admits an asymptotically cylindrical
Ricci-flat K\"ahler metric with holonomy~$SU(3)$.
Note that the cylindrical asymptotic model for this metric is determined by
the Ricci-flat K\"ahler structure in the K\"ahler class of $D$ in $X_8$.

We would like to choose two octic threefolds $X_8^{(i)}$, $i=1,2$ and a K3
surfaces $D_i$ in each, so as to satisfy the Donaldson matching condition.
We do this by applying some general theory of K3 surfaces and their moduli
(see \cite[Ch. VIII]{BHPV}). The key point is that one can determine a
Ricci-flat K\"ahler K3 surface $D$, up to isomorphism, by a data of the
integral second cohomology $H^2(D,\ZE)$.

Recall that all K3 surfaces are
diffeomorphic and the intersection form makes $H^2(D,\ZE)$ into a lattice.
There is an isomorphism, called a marking, $p:H^2(D,\ZE)\to L$ to a fixed
non-degenerate even unimodular lattice $L$ with signature $(3,19)$. We shall
refer to $L$ as the {\em K3 lattice}; its bilinear form is given by the
orthogonal direct sum $L=3H\oplus 2(-E_8)$ of 3 copies of the hyperbolic plane
lattice $H=\bigl(\begin{smallmatrix} 0&1\\ 1&0 \end{smallmatrix}\bigr)$ and
2 copies of the negative definite root lattice $-E_8$ of rank~8. 
Now if $D\subset\CP^5$ is an octic K3 surface, then the image 
$p(\kappa_I)$ of the K\"ahler class of $D$ is primitive (non-divisible in~$L$
by an integer $>1$) and $p(\kappa_I)\cdot p(\kappa_I)=8$, computed in the
bilinear form of $L$. The images $p(\kappa_J),\, p(\kappa_K)$ span a positive
2-plane $P$ orthogonal to $p(\kappa_I)$ in the real vector space $L\otimes\RE$.
Conversely, the positive 2-planes $P$ arising in this way form a dense open
set in the Grassmannian of positive 2-planes orthogonal to $p(\kappa_I)$ in
$L\otimes\RE$.

It is known that the group of lattice isometries of $L$ acts
transitively on the set of all primitive vectors with a fixed value of
$v\cdot v$ in~$L$. We can therefore choose two octic K3 surfaces with
hyper-K\"ahler structures (in the respective K\"ahler classes)
$(D_1;\kappa'_I,\kappa'_J,\kappa'_K)$,
$(D_2;\kappa''_I,\kappa''_J,\kappa''_K)$ and the markings $p_1$, $p_2$ with
$p_1(\kappa'_I)=p_2(\kappa''_J)$, $p_1(\kappa'_J)=p_2(\kappa''_I)$ in $L$,
and $p_1(\kappa'_K)=-p_2(\kappa''_K)$ in $L\otimes\RE$ thereby achieving a
matching.

Choosing the ambient octic threefolds $X'_8,X''_8$ for the latter $D_1,D_2$,
blowing up these threefolds to obtain asymptotically cylindrical Ricci-flat
threefolds by Theorem~\ref{bCY}, and applying Theorem~\ref{gluing} to the
respective connected sum, we obtain a simply-connected compact 7-manifold $M$
with a metric of holonomy~$G_2$.

We may consider in a very similar way, in place of one of both $X_8$'s above, a
smooth intersection $X_6$ of a quadric and a cubic in $\CP^5$.  The respective
K3 divisor then is an intersection of a quadric and a cubic in $\CP^4$
and the image of the K\"ahler class of this divisor has square 6 in the
bilinear form~$L$.
\end{example}

More generally, it was shown in~\cite[\S 6,7]{g2paper} that in place of $X_8$,
$X_6$ in the above example we can consider any non-singular {\em Fano
threefold}~$V$, i.e.\ a projective complex 3-dimensional manifold such that
the image of the first Chern class $c_1(V)$ in the de Rham cohomology can be
represented by some K\"ahler form on~$V$. Equivalently, the anticanonical
bundle $K^{-1}_V$ is ample. Smooth Fano threefolds are completely classified;
up to deformations, there are 105 algebraic types \cite{Isk,MM}.

Every Fano threefold $V$ is simply-connected and a
generic anticanonical divisor $D$ in~$V$ is a (smooth) K3 surface~\cite{sh}.
A threefold $\bW$ is obtained by blowing up a connected complex curve
representing the self-intersection cycle $D\cdot D$ (in the sense of the Chow
ring). Then $\bW$ and the proper transform of $D$ satisfy the hypotheses of
Theorem~\ref{bCY}. Alternatively, if $D\cdot D$ is represented by a finite
sequence of curves, then $\bW$ may be defined by successively blowing up these
curves. We shall refer to any such threefold $\oW$ to be of {\em Fano type}.

 A K\"ahler K3 surface $D$ and its proper transform in $\bW$ can be
assumed isomorphic by choosing an appropriate K\"ahler metric on~$\bW$. Then
the cylindrical asymptotic model for~$W$ is determined by the K3 surface~$D$
with the K\"ahler metric restricted from $V$.

For a general Fano $V$, the class of anticanonical K3 surfaces $D$ arising in
the deformations of~$V$ will correspond to an open dense subset of
{\em lattice-polarized} K3 surfaces. This latter class is defined by the
condition that the Picard lattice $H^{1,1}(D,\RE)\cap H^2(D,\ZE)$ contains a
sublattice isomorphic to a fixed lattice $N$ and this sublattice contains a
class of some K\"ahler form. In the case of algebraic K3 surfaces 
of a fixed degree, as in the example above, $N$ is generated by the K\"ahler
form $\kappa_I$ induced from the embedding of $D$ in the projective space. In
general, $N$ arises as $\iota^* H^2(V,\ZE)$ from the embedding $\iota:D\to V$.
The rank of $N$ is the Betti number $b^2(V)$ as $\iota^*$ in injective by the
Lefschetz hyperplane theorem.
\smallskip

Another source of examples for Theorem~\ref{bCY} was given by Lee and the
author in~\cite{k3g2}. The construction uses K3 surfaces $S$ with
{\em non-symplectic involution}, a holomorphic map $\rho:S\to S$, such that
$\rho^*$ restricts as $-1$ on $H^{2,0}(S)$.  The K3 surfaces of this type were
completely classified up to deformation by Nikulin~\cite{nikulin2}, who
determined the complete system of invariants and fixed point set of $\rho$ for
each deformation family. We require the fixed point set of $\rho$ to be
non-empty; this occurs in all but one of the 75 deformation families.

Let $\psi:\CP^1\to\CP^1$ denote the holomorphic involution
$\psi(z_0:z_1)=z_1:z_0$ fixing exactly two points. The quotient
$Z=(S\times\CP^1)/(\rho,\psi)$ is then an orbifold whose singular locus is a
disjoint union of smooth curves. The desired 3-fold $\bW$ is defined by the
resolution of singularities diagram for~$Z$,
$$
\begin{CD}
\widetilde{W} 	@>>> \oW\\
@VVV 		@VVV\\
S\times\CP^1 	@>>> Z,
\end{CD}
$$
where the vertical arrows correspond to blowing up the fixed locus of
$(\rho,\psi)$ in $S\times\CP^1$ and the singular locus of~$Z$  and the
horizontal arrows are the quotient maps.

The anticanonical divisor $D$ in~$\bW$ arises as the (pre-)image of
$S\times\{p\}$, via the above diagram, where $\psi(p)\neq p$. Such $D$ is
clearly isomorphic to the K\"ahler K3 surface~$S$ and has trivial normal
bundle in~$\bW$. It can be checked that $\bW$ and $W=\bW\setminus D$ are
simply-connected (the condition that $\rho$ have fixed points is needed
here). Thus $W$ has an asymptotically cylindrical Ricci-flat K\"ahler metric
by Theorem~\ref{bCY}.

The pull-back $\iota^*:H^2(\bW,\ZE)\to H^3(D,\ZE)$ defined by the embedding of
$D$ makes $D$ into a lattice polarized K3 surface with $N$ corresponding to
the sublattice of all classes fixed by $\rho^*$ in $H^2(D,\ZE)$.
On the other hand, $\iota^*$ has a kernel of dimension at least~2. A
threefold $\bW$ obtained from K3 surface with non-symplectic involution is
therefore {\em never} deformation equivalent to any threefold of Fano type
(assuming $D\cdot D$ in the latter threefold was represented by a single
curve).
\smallskip

The matching problem in all the examples becomes entirely a consideration on
the K3 lattice~$L$, as illustrated by the example in the beginning of this
subsection. In general, the argument is more technical and requires results on
the lattice embeddings~\cite{nikulin1}.

One simple sufficient (though not necessary) condition for the existence of
the Donaldson matching for representatives in the two classes of lattice
polarized K3 surfaces is that the rank of each polarizing lattice $N_i$ is
$\le 5$.
\smallskip

All the irreducible $G_2$-manifolds $M$ constructed from threefolds in the
above examples are simply-connected.
The cohomology of compact irreducible $G_2$-manifolds $M$ coming from the
connected sum construction may be determined by application of the
Mayer--Vietoris exact sequence and generally depends on the choice of
matching. However, the sum of the Betti numbers
\begin{equation}\label{b2b3}
b^2(M)+b^3(M)=b^3(\oW_1)+b^3(\oW_2)+2d_1+2d_2+23,
\end{equation}
for any matching, depends only on the threefolds $\oW_i$ and the dimensions
$d_i$ of the kernel of $\iota^*:H^2(\bW,\RE)\to H^3(D,\RE)$.
The quantities in~\eqref{b2b3} can be determined by standard methods
(adjunction formula, Lefschetz--Bott hyperplane theorem) from known
algebraic invariants of Fano threefolds or, respectively, of non-symplectic
involutions.

In particular, the Fano threefold $X_8$ discussed in Example~\ref{one} above
has $b^2(X_8)=1$, $b^3(X_8)=28$ and its blow-up has $b^2(\bW)=2$,
$b^3(\bW)=38$. An irreducible compact $G_2$-manifold $M$ constructed from a
pair of $X_8$'s has $b^2(M)=0$ and then $b^3(M)=99$ as $d_i$ vanish in this
case. This irreducible $G_2$-manifold is topologically distinct from
the examples given by Joyce via resolution of singularities;
the only irreducible $G_2$-manifold in~\cite{joyce-book} with $b^2=0$ has
$b_3=215$.
The property $b^2(M)=0$ holds in many other examples coming from pairs
of threefolds of Fano type and these latter $G_2$-manifolds typically have
smaller $b^2$ and larger $b^3$ than the examples given by Joyce.
(Note also that every compact irreducible $G_2$-manifold $M$ must have
$b^1(M)=0$ by Proposition~\ref{irreducible} but $b^3(M)$ cannot vanish as the
$G_2$ 3-form $\phi$ on $M$ is harmonic.)
\medskip

Corti, Haskins, Nordstr\"om and Pacini~\cite{CHNP1,CHNP2} generalized the
class of threefolds of Fano type by considering {\em weak Fano} threefolds $V$
whose anticanonical bundle $K^{-1}_V$ is only required to be big and nef.
(Every such $V$ may be obtained as resolution of an
appropriate singular Fano threefold.) They identified a large subclass called
semi-Fano threefolds and generalized for this class the properties required
in the construction of $G_2$-manifolds from threefolds $W$ of Fano type.
This generalization dramatically increased the number of examples of connected
sum $G_2$-manifolds. Some of the examples were shown to be 2-connected which
allows to determine their diffeomorphism type by computing certain standard
invariants.

More recently, Braun \cite{braun} gave a toric geometry construction, from
certain lattice polytops, of examples of pairs $\bW,D$ defining asymptotically
cylindrical Calabi--Yau threefolds by Theorem~\ref{bCY}. Useful invariants of
$\bW$ e.g.\ the Hodge numbers can be computed by combinatorial formulae.

Nordstr\"om \cite{extra-twisted} gave an interesting generalization of the
connected sum construction, by replacing~\eqref{rotation} with a different
`hyper-K\"ahler rotation' and taking finite quotients of asymptotically
cylindrical Calabi-Yau threefolds $W$.
Applications of the construction include topologically new examples of compact
irreducible $G_2$-manifolds some of which have a non-trivial finite
fundamental group. 
\smallskip

In conclusion, we mention two works which contain results concerning 
relations between the two types of construction of $G_2$-manifolds discussed
in these notes.

Nordstr\"om and the author identified in \cite{KN} an example of a compact
irreducible $G_2$-manifold given by Joyce~\cite{joyce-book} where the
underlying 7-manifold is diffeomorphic to one obtainable from the construction
in~\cite{k3g2}. Further, the two respective families of $G_2$-metrics on this
manifold are connected in the $G_2$-moduli space.

On the other hand, some of the $G_2$-manifolds given by Joyce cannot possibly
be obtained by the connected sum construction. The result is due to Crowley
and Nordstrom~\cite{CN} who constructed an invariant of $G_2$-structures which
is equal to 24 for each connected sum~\eqref{7mfd} but is odd for some
examples in~\cite{joyce-book}.

\end{document}